\documentstyle[12pt,righttag]{amsart}
\numberwithin{equation}{subsection}
\theoremstyle{plain}
  \newtheorem{th}{Theorem}[subsection]
  \newtheorem{pr}[th]{Proposition}
  \newtheorem{lem}[th]{Lemma}
  \newtheorem{cor}[th]{Corollary}
\theoremstyle{definition}
  \newtheorem{df}[th]{Definition}
\theoremstyle{remark}
 \newtheorem{rem}[th]{Remark}

\def\b{{\frak b}}

\def\e{{\epsilon}}
\def\g{{\frak{gl}}}
\def\h{{\frak t}}

\def\l{{\ell}}

\def\n{{\frak n}}

\def\t{{{\frak t}}}

\def\ob{{obj\,}}

\def\Vn{{{V^{\*\l}_n}}}

\def\C{{\Bbb C}}

\def\K{{\cal M}}
\def\Sym{{\frak S}}
\def\L{{\cal L}}

\def\O{{\cal O}}

\def\W{{W}}
\def\Y{{\cal Y}}
\def\Z{{\Bbb Z}}
\def\+{\mathop{\oplus}}
\def\*{\mathop{\otimes}}

\def\pair{{\lm,\mu}}
\def\one{{\text{\rom{\bf 1}}}}

\def\gl{{{\frak{gl}}}}
\def\lm{\lambda}
\def\al{\alpha}
\def\Rep{{{\cal R}}}
\def\ch{^\vee}

\def\bra{{\langle}}
\def\ket{{\rangle}}
\def\End{{\text{\rom{End}}}\,}
\def\Hom{{\text{\rom{Hom}}}\,}
\def\dim{{\text{\rom{dim}}}\,}

\def\gen{{\text{\rom{gen}}}}
\def\Ker{{\text{\rom{Ker}}}}

\def\gr{{\text{\rom{gr}}}}
\def\tr{{{\text{\rom{tr}}}}}
\def\invo{{\iota}}
\begin{document}
\title[degenerate affine Hecke algebras]
{Rogawski's conjecture on the Jantzen filtration 
for
the degenerate affine Hecke algebra 
of \\ type $A$} 
\author{Takeshi Suzuki}
\address{Research Institute for Mathematical Sciences, 
Kyoto University, Japan}
\email{takeshi@@kurims.kyoto-u.ac.jp}
\thanks{The author is supported by the JSPS 
Research Fellowships for Young Scientists.}
\maketitle
\begin{abstract}
The functors constructed by Arakawa and the author
relate the representation theory of
$\g_n$ and that of the
degenerate affine Hecke algebra $H_\l$ of $GL_\l$.
They transform the Verma modules over $\g_n$
to the standard modules over $H_\l$. 
In this paper we prove that
they transform the simple modules 
to the simple modules (in more general situations 
than in the previous paper). 
We also prove that they 
transform the Jantzen filtration on the Verma modules
to that on the standard modules.
We obtain the following results
for the representations of $H_\l$
by translating the corresponding results 
for $\g_n$ through the functors:
(i) the (generalized) Bernstein-Gelfand-Gelfand
 resolution for a certain class of
simple modules,
(ii) the multiplicity formula for 
the composition series of the standard modules,
and (iii) its refinement concerning the Jantzen filtration
on the standard modules,
which was conjectured by Rogawski.
\end{abstract}
\begin{center}
  {\sc Introduction}
\end{center}
This paper is a continuation of 
the paper\ \cite{AS}, in which 
we gave
functors from $\O(\g_n)$ to $\Rep(H_\l)$.
Here $\O(\g_n)$ denotes the Bernstein-Gelfand-Gelfand
(in short, BGG) category of representations of
the complex Lie algebra
$\g_{n}$, and
$\Rep(H_\l)$ denotes the category of finite-dimensional 
representations of 
the degenerate affine Hecke algebra $H_\l$ of $GL_\l$
introduced by Drinfeld\ \cite{Dr}.

Let us review the results in \cite{AS}.
Let $\h_n^*$ and $W_n$ denote the space of weights and Weyl group
of $\g_n$ respectively.
For $\lm\in \h^*_n$,
let $M(\lm)$ denote the Verma module with highest weight $\lm$
and $L(\lm)$ its simple quotient.
Let $V_{n}=\C^n$ denote the vector representation of $\g_{n}$. 
For each $\lm\in\h^*_n$ and
$X\in\O(\g_{n})$, 
we define  an action of $H_\l$
on the finite-dimensional vector space
$F_\lm(X)=\Hom_{\g_{n}}(M(\lm),X\*\Vn)$. 
Under the condition that
$\lm+\rho$ is dominant, 
we proved that
the functor $F_\lm$ is exact and
$F_\lm(M(\mu))$ is isomorphic to
$\K(\lm,\mu)$ unless it is zero.
Here $\K(\pair)\in \Rep(H_\l)$ denotes
the standard module.
With the restriction  $\l=n$, we proved that 
$F_\lm(L(\mu))$ is isomorphic to 
the unique simple quotient $\L(\pair)$ 
of $\K(\pair)$ unless it is  zero.
Any simple $H_\l$-module is thus obtained.
To prove the irreducibility
of $F_\lm(L(\mu))$, we compared the multiplicities
of the simple modules in the composition series of
$M(\mu)$ and those in $\K(\pair)$
by using
the Kazhdan-Lusztig type multiplicity formulas
known for $\O(\g_{n})$ 
and $\Rep(H_\l)$. (See (b)~(c) below.) 

In the present paper,  further properties of the functors 
are deduced from the
key observation 
that 
the $\g_{n}$-contravariant bilinear form on a highest weight
$\g_n$-module $X$
induces the $H_\l$-contravariant bilinear form on $F_\lm(X)$. 
The irreducibility of $F_\lm(L(\mu))$ is deduced from
the non-degeneracy of the bilinear form. 
As a consequence, we can determine the images of simple 
$\g_n$-modules (Theorem \ref{th;image_of_simple})
 without assuming  $\l=n$
or referring to the multiplicity formulas.

We also prove that $F_\lm$ transforms
the Jantzen filtration on $M(\mu)$ to
that on $F_\lm(M(\mu))\cong \K(\pair)$ (Theorem \ref{th;Jantzen}).

The followings are the consequences of these results.

\smallskip\noindent
(i) We obtain a resolution for a certain class 
of simple $H_\l$-modules by
applying $F_\lm$ to the BGG resolution \cite{BGG}
and its generalization by Gabber-Joseph \cite{GJ;BGG}
for $\g_n$-modules.
This generalizes the results of
Cherednik \cite{Ch;character} and Zelevinsky \cite{Zel;resolvant}.

\medskip\noindent
(ii) To simplify the descriptions, we assume
 $\lm$ and $\mu$ are dominant integral weights.
(More general cases are treated in \S \ref{ss;KL}.)
Set $w\circ\mu=w(\mu+\rho)-\rho$ for $w\in W_n$ and
let $w,y\in W_n$ be such that $\lm-w\circ\mu$ and 
$\lm-y\circ\mu$ are weights of $\Vn$.
We have a direct proof 
of the following formula:
\begin{equation*}\label{eq;mult}
[M(w\circ\mu):L(y\circ\mu)]=
[\K(\lm,w\circ\mu):\L(\lm,y\circ\mu)].\tag{a}
\end{equation*}

Let $P_{w,y}(q)$ denote the Kazhdan-Lusztig polynomial of
$W_n$.
The formula (a) implies the equivalence of the
following two multiplicity formulas:
\begin{alignat}{2}
  [\K(\lm,w\circ\mu)&:\L(\lm,y\circ\mu)]&=&P_{w,y}(1),\tag{b}\\
   [M(w\circ\mu)&:L(y\circ\mu)]&=&P_{w,y}(1).\tag{c}
\end{alignat}
The formula (b) was 
proved by Ginzburg\ \cite{Gi} (see also \cite{CG})
for affine Hecke algebras,
and (c) was
proved by Beilinson-Bernstein\ \cite{BB}
and Brylinski-Kashiwara\ \cite{BK} by 
using the geometric method and the theory of perverse sheaves.
We remark that our proof of (a) is purely algebraic.

\medskip\noindent
(iii) We have a refinement of the formula (a):
Let $\lm,\mu$ and $w,y$ be as in (ii).
(See \S \ref{ss;Rog} for more general cases.)
Let 
\begin{align*}
M(\mu)&=M(\mu)_0\supseteq M(\mu)_1 \supseteq
M(\mu)_2\supseteq
 \cdots,\\
\K(\lm,\mu)&=\K(\pair)_0\supseteq \K(\pair)_1 \supseteq \K(\pair)_2 \supseteq 
\cdots
\end{align*}
be the Jantzen filtrations on $M(\mu)$ and $\K(\pair)$, 
respectively.
Since $F_\lm$ preserves the Jantzen filtration, we have
\begin{equation*}\label{eq;mult'}
  [M(w\circ\mu)_j:L(y\circ\mu)]=
[\K(\lm,w\circ\mu)_j:\L(\lm,y\circ\mu)].\tag{a'}
\end{equation*}

The Jantzen filtration on 
standard modules
over affine Hecke algebras of $GL$
was introduced by Rogawski~\cite{Ro}.
He conjectured
a refinement of the formula (b)
concerning the Jantzen filtration.
Rogawski's conjecture 
was proved by Ginzburg\footnote{
I.~Grojnowsky 
announced similar results
in a series of his lectures at Kyoto 1997. 
He also treated affine Hecke algebras at a root of unity.}.
(The result is announced
in \cite{Gi;ICM} without details.)
A degenerate affine Hecke analogue of
Rogawski's conjecture is written as follows:
\begin{equation*}
  \sum_{i \in \Z_{\geq 0}} 
[\gr_i \K(\lm,w\circ\mu):\L(\lm,y\circ\mu)]
q^{(l(y)-l(w)-i)/2}=P_{w,y}(q).\tag{b'}
\end{equation*}
The formula (a') implies the equivalence between
(b') and 
the improved Kazhdan-Lusztig multiplicity formula
\begin{equation*}
\sum_{i \in \Z_{\geq 0}} 
[\gr_i M(w\circ\mu):L(y\circ\mu)]q^{(l(y)-l(w)-i)/2}
=P_{w,y}(q), \tag{c'}
\end{equation*}
which was proved in \cite{BB;Jantzen}.

\medskip\noindent
{\bf Acknowledgment.}
I would like to thank T.\ Arakawa and A.\ Tsuchiya for
the collaboration in \cite{AST} and \cite{AS},
which lead to this work.
I would like to thank I.\ Grojnowski for
his beautiful lectures on affine Hecke algebras,
which 
motivated the author for the present problems. 
Thanks are also due to K.\ Iohara, T.\ Miwa,
Y.\ Koyama and Y.\ Saito
for valuable discussions.
I am grateful to M.\ Kashiwara for careful reading of 
the manuscript and a lot of important comments.
\section{Basic definitions}
\subsection{Lie algebra $\g_n$}\label{ss;sl}
Let $\gl_n$ denote
the Lie algebra consisting of
all $n\times n$ matrices with entries in $\C$.
Let $\t_n$ be the Cartan subalgebra of $\gl_n$
consisting of all diagonal matrices.
An inner product is defined on $\gl_n$ by
\begin{equation}\label{eq;inner}
  (x|y)_n=\tr(xy)
\end{equation}
for $x,y\in\gl_n$.
Let $\t_n^*$ denote
the dual space of $\t_n$. 
The natural pairing is denoted by
$\bra~,~\ket_n : \t_n^*\times\t_n\to\C$. 
Let $E_{i,j}$ $(1\leq i,j \leq n)$
denote the matrix with only nonzero entries $1$
at the $(i,j)$-th component.
Define 
a basis $\{\e_i\}_{i=1,\dots,n}$ of $\t_n^*$
by $\e_i(E_{j,j})=\delta_{i,j}$, and
define the roots by $\al_{ij}=\e_i-\e_j$ and
the simple roots by
$\al_i=\e_i-\e_{i+1}$.

Put 
\begin{align}
  R_{n}&=\{\al_{ij}\mid  1\leq i\neq j\leq {n}\},\\
  R_{n}^+&=\{\al_{ij}\mid  1\leq i<j\leq {n}\},
\quad R_{n}^-=R_{n}\setminus R_{n}^+,\\
  \Pi_{n}&=\{\al_i\mid i=1,\dots{n}-1\}.
\end{align}
Then $R_{n}\subseteq \t^*_{n}$ 
is a root system of type $A_{{n}-1}$.
Since the restriction of  $(~|~)_n$ to $\t_n$
is non-degenerate,
we have an isomorphism $\t_n^*
\stackrel{\sim}{\to}
\t_n$, 
whose image of $\xi\in\t_n^*$ is denoted by $\xi\ch$.
In particular we have
 $\e\ch_i=E_{i,i}$ and $\al\ch_i=E_{i,i}-E_{i+1,i+1}$.

Putting $\n_n^+=\+_{i<j}\C E_{i,j}$, $\n_n^-=\+_{i>j}\C E_{i,j}$,
we have a triangular decomposition
 $\g_{n}=\n^+_n\+\h_{n}\+\n^-_n$ .
We put $\b_n^\pm=\n_n^\pm\+\h_n$.

Let $\sigma$ denote the involution on $\gl_n$ given by
the transposition: $\sigma(E_{i,j})=E_{j,i}$.
The inner product $(~|~)_n$ is invariant with respect to $\sigma$:
$(\sigma(x)|\sigma(y))_n=(x|y)_n$ for all $x,y\in\g_n$.

Put $\rho=\frac{1}{2}\sum_{\al\in R_n^+}\al$ and define
\begin{align}
 Q_{n}&=\+_{i=1}^{{n}-1}\Z\al_i,\\
 D_{n}&=\{ \lm\in \t_{n}^*\mid
\bra \lm+\rho ,\al\ket_n \notin \Z_{< 0}
\text{ for all }\al\in R_{n}^+\}, \\
  D_{n}^\circ&=\{ \lm\in \h_{n}^*\mid
\bra \lm , \al\ket_n \notin \Z_{< 0}
\text{ for all }\al\in R_{n}^+\}, \label{eq;dominant_weight}\\
 P_{n}&=\+_{i=1}^n\Z \e_i,\quad  P_{n}^+=P_{n}\cap D_{n}^\circ.
\end{align}
An element of $D_{n}^{\circ}$ (resp. $P_{n},\, P^+_{n}$)
is called a {\it dominant} 
(resp.{\it integral, dominant integral})
weight.
\subsection{Weyl group}\label{ss;Weyl_group}
Let $W_{n}\subset \rom{GL}(\t^*_{n})$ be the Weyl group associated to
the root system $(R_{n}, \Pi_{n})$, 
which is by definition generated by
the reflections $s_{\al}$ $(\al\in R_{n})$ 
defined by 
\begin{equation}
  s_\al(\lm) =\lm-  \bra\lm, \al\ch\ket_n \al \quad (\lm\in\t^*_{n}).
\end{equation}
We often write
$s_i=s_{\al_i}$ for $\al_i\in \Pi_{n}$.
Note that $W_{n}$
is isomorphic to the symmetric group $\Sym_{n}$.

We often use another action of $W_{n}$ on $\h^*_{n}$,
which is given by
\begin{equation}
  \label{eq;dot_action}
  w\circ\lm=w(\lm+\rho)-\rho\quad (w\in W_{n},\ \lm\in\h_{n}^*).
\end{equation}
For $w,y\in W_{n}$, we write $ w\geq y$
if and only if  $y$ can be obtained as a subexpression of a 
reduced expression of $w$.
The resulting relation in $W_{n}$
defines a partial order called the {\it Bruhat order}.
\subsection{Representations of $\g_n$}
%
For a $\h_{n}$-module $X$ and $\lm\in\h^*_{n}$, put
\begin{align}
 X_{\lm}         &=\{ v\in X\mid hv=\bra \lm,h\ket_n v 
\text{ for all }h\in\h_{n}\},\label{eq;weight_space}\\
 X_{\lm}^\gen    &=\{ v\in X\mid  (h-\bra\lm,h\ket_n)^kv=0
\text{ for all }h\in \t_n,
\text{ some } k\in\Z_{>0}
\},\label{eq;generalized_weight_space}\\
 P(X)          &=\{\lm\in\h^*_{n}\mid X_{\lm}\neq 0\}.
\label{eq;set_of_weight}
\end{align}
The space $X_{\lm}$ (resp $X_{\lm}^\gen$) 
is called the {\it weight space} (resp. {\it generalized weight space})
of weight $\lm$
with respect to $\h_{n}$, and 
an element of $P(X)$ is called a weight of $X$.

Let $U(\g_n)$ denote the universal enveloping algebra of $\g_n$.
Let 
$\O=\O(\g_{n})$ denote
the category of $\g_n$-modules 
which are finitely generated over $U(\g_n)$,
$\n^+_n$-locally finite 
and $\h_{n}$-semisimple (see \cite{BGG}).
The category $\O$ is closed under the operations such as
forming subquotient modules,
finite direct sums, and tensor products with finite-dimensional modules.
For $\lm\in\h^*_{n}$,
let $M(\lm)=U(\g_{n})\*_{U(\b^+_{n})}\C v_{\lm}$ 
denote 
the Verma module with highest weight $\lm$, where 
$v_\lm$ denotes the
highest weight vector.
The unique simple quotient of $M(\lm)$ is denoted by $L(\lm)$.
The modules $M(\lm)$ and $L(\lm)$ are objects of $\O$.

Let $\chi_{\lm}:Z(U(\g_{n}))\to \C$ denote
the infinitesimal character of $M(\lm)$.
We introduce an equivalence relation in $\h^*_n$ by
\begin{equation}
\lm\sim\mu\Leftrightarrow \lm=w\circ\mu\ \text{ for some }
w\in W_{n}.\label{eq;equivalence_relation}
\end{equation}
Then it follows that
$\chi_\lm= \chi_\mu$ if and only if $\lm\sim\mu$.
Let $[\lm]$ denote the equivalence class of $\lm\in\h_n^*$.
Define the full subcategory $\O_{[\lm]}$ of $\O$
by
$$\ob \O_{[\lm]}=\{X\in 
\ob \O\mid (\Ker\chi_\lm)^k X=0 \text{ for some } k\}.$$ 
Then any $X\in\ob\O$ admits a decomposition
\begin{equation}
  \label{eq;decomposition}
X=\+_{[\lm]\in \h_n^*/ \sim }X^{[\lm]}  
\end{equation}
such that $X^{[\lm]}\in \ob\O_{[\lm]}$.
The correspondence $X\mapsto X^{[\lm]}$
gives an exact functor on $\O$.
\begin{lem}  \label{lem;bij}
Let $\lm\in D_{n}$. 
Then the natural map 
$ 
(X^{[\lm]})_\lm
\to(X/\n^-_nX)_{\lm}
$ 
is bijective.
\end{lem}
\begin{rem}
\rom{(i)} 
There also exists a canonical bijection
$\Hom_{\g_n}(M(\lm),X)\cong 
(X^{[\lm]})_{\lm}$
for $\lm\in D_n$.

\smallskip\noindent
\rom{(ii)}
A proof of Lemma \ref{lem;bij}
for integral $\lm$ is given in \cite{AS}.
The generalization to non-integral cases is similarly proved.
\end{rem}
\section{Degenerate affine Hecke algebras and their representations}
\subsection{Degenerate affine Hecke algebras}\label{ss;DAHA}
For a group $G$, let  $\C[G]$
denote its group ring.
Let $S(\t_\l)$ denote the symmetric algebra of $\t_{\l}$, 
which is isomorphic to
the polynomial ring $\C[\e_1\ch,\dots,\e_\l\ch]$.
\begin{df}
The {\it degenerate $($or graded$)$ affine Hecke algebra} $H_{\l}$ 
of $GL_{\l}$ is the unital associative algebra over $\C$ defined by
the following properties:

\medskip\noindent
(i) As a vector space, $H_{\l}\cong \C[W_{\l}]\* S(\t_\l)$. 

\smallskip\noindent
(ii) The subspaces $\C[W_{\l}]\*\C$ and $\C\*S(\t_\l)$ are
subalgebras of $H_{\l}$ in a natural fashion (their images 
will be identified with $\C[W_{\l}]$ and $S(\t_\l)$ respectively).

\smallskip\noindent
(iii) The following relations hold in $H_{\l}$: 
\begin{equation}\label{eq;relation}
  s_\al\cdot \xi-s_\al(\xi)\cdot s_\al=-\bra \al,\xi\ket_\l\quad
(\al\in \Pi_{\l},\ \xi\in\t_{\l}).
\end{equation}
\end{df}
It is easy to verify the following lemma.
\begin{lem}\label{lem;involution}
There exists a unique anti-involution $\invo$ on $H_\l$
such that
$$
\invo(w)=w^{-1}\ (w\in W_\l),\ \ 
\invo(\xi)=\xi\ (\xi\in\t_\l).$$
\end{lem}
For a subset $B\subseteq\Pi_\l$, 
let $\t_B$ denote the subspace of $\t_\l$ spanned by 
all $\e_i\ch$ such that $\bra \al,\e_i\ch\ket\neq 0$
for some $\al\in B$.
Put
\begin{equation}\label{eq;parabolic_subalgebra}
H_B=\C[W_B]\* S(\t_B)\subseteq H_\l.
\end{equation}
Then it turns out that $H_B$ is a subalgebra of $H_\l$.
 %

%
\subsection{Induced modules}\label{ss;induced}
For a pair $\Delta=[a,b]$ 
of complex numbers 
 such that $b-a+1=\l\in\Z_{\geq0}$, there exists
 a unique one-dimensional representation 
$\C_{\Delta}=\C \one_{\Delta}$ of $H_{\l}$ (we put $H_0=\C$ for
convenience)
such that 
\begin{align}
  \label{eq;condition1}
  w \one_{\Delta} &=\one_{\Delta}\quad (w\in W_{\l}), \\
  \label{eq;condition2}
  \e\ch_i\one_{\Delta}&= (a+i-1)\one_{\Delta}\quad (i=1,\dots,\l).
\end{align}

Let $\lm\in\h_n^*$
and
$\mu\in\h_n^*$
be such that $\lm-\mu\in P(\Vn)$.
Then putting 
\begin{align}\label{eq;partition}
\l_i&=\bra\lm-\mu,\e\ch_i\ket_n\in \Z_{\geq 0}\quad (i=1,\dots,n),
\end{align}
 we have
$ \l= \sum_{i=1}^{n}\l_i$.
Let $\Pi_{\lm,\mu}\subseteq\Pi_\l$ be
the subset associated to the partition $(\l_1,\dots,\l_n)$:
$\Pi_{\lm,\mu}=\{\al_i\mid i\neq \sum_{k=1}^j \l_k\ \text{for any }j\}$.  
Put 
\begin{align}
W_{\lm,\mu}&=W_{\Pi_{\lm,\mu}}\subseteq W_\l,\quad
H_{\lm,\mu}=H_{\Pi_{\lm,\mu}}\subseteq H_\l.
\end{align}
Note that $H_{\lm,\mu}= H_{\l_1}\*\cdots\* H_{\l_n}
=S(\t_\l)\*\C[W_{\lm,\mu}]$.
Put
\begin{align}
  \Delta_i&=[\bra\mu+\rho,~
\e\ch_i\ket_n, \bra\lm+\rho,\e_i\ch\ket_n-1]
\in\C^2.
\end{align}
Define the parabolically induced module 
$\K(\pair)$ associated to $(\pair)$ by
\begin{equation}
\label{eq;parabo}
  \K(\pair)=H_\l\*_{H_{\lm,\mu}}
(\C_{\Delta_1}\*\cdots\* \C_{\Delta_{n}}).
\end{equation}
Evidently $\K(\pair)$ is a cyclic module with a cyclic weight vector 
\begin{equation}\label{eq;one}
  \one_{\pair}:=
  \one_{\Delta_1}\*\cdots\* \one_{\Delta_k},
\end{equation}
whose weight $\zeta_\pair$ is given by
\begin{equation}\label{eq;weight_of_one}
  \bra\zeta_\pair,\e\ch_j\ket_\l=
\bra \mu+\rho,\e_i\ch\ket_n+j-\sum_{k=1}^{i-1}\l_k-1
\quad
\text{for}\quad\sum_{k=1}^{i-1}\l_k  < j\leq\sum_{k=1}^{i}\l_k.
\end{equation}
It is also obvious that 
$\K(\pair)\cong \C[W_{\l}/ W_{\lm,\mu}]$
as a $\C[\W_{\l}]$-module and thus its dimension is given by
$\dim \K(\pair)=
{\l !}/({\l_{1}!\cdots \l_{k}!}).$
Recall that the simple modules of $W_\l$
are parameterized by unordered partitions of $\l$ (or Young diagrams of size 
$\l$).
We let $S_\gamma$ denote the simple $W_\l$-module
corresponding to the partition $\gamma$.
Let $[\lm-\mu]$ 
denote the unordered partition of $\l$
obtained from $(\l_1,\dots,\l_k)$ by forgetting the order.
As is well-known, it holds that 
\begin{equation}\label{eq;decomposition_as_W-module}
  \K(\pair)\cong S_{[\lm-\mu]}\oplus
\bigoplus_{\beta\triangleright [\lm-\mu]}S_\beta^{\+ a_\beta},
\end{equation}
as a $\C[W_\l]$-module .
Here
$\triangleright$ denotes the dominance order
in the set of partitions, 
and $a_\beta$ are some non-negative integers.

Let $\Y_\l(n)$ denote the set of Young diagrams of size $\l$ 
consisting of at most $n$ rows.
We say that an $H_\l$-module $Y$ is of {\it level} $n$
if 
$Y=\+_{\gamma\in \Y_\l(n)}S_{\gamma}^{\+ a_\gamma}$ 
for some $a_\gamma\in \Z_{\geq0}$.
The induced module $\K(\pair)$ $(\lm,\mu\in\h_n^*)$ is
of level $n$.
Of course, any finite-dimensional $H_\l$-module is of level $\l$.
\subsection{Zelevinsky's classification of simple modules}
The representation theory
 of the degenerate affine Hecke algebra
is related to that of the affine Hecke algebra
by Lusztig \cite{Lu}.
Thus the statements in this subsection are deduced from 
\cite[Theorem 6.1]{Zel;induced}
and \cite[\S 5]{Ro}.
(See also \cite{Ch;special_bases}.)

\begin{th}[\cite{Zel;induced,Ro}]
\label{th;usq}
Let $\lm\in D_{n}$ and $\mu\in \lm-P(\Vn)$.

\medskip\noindent
\rom{(i)}
In the decomposition \eqref{eq;decomposition_as_W-module},
$S_{[\lm-\mu]}$ generates $\K(\pair)$ over $H_\l$.

\smallskip\noindent
\rom{(ii)}
The $H_{\l}$-module $\K(\pair)$ has the unique simple quotient, 
which is
denoted by $\L(\pair)$.

\smallskip\noindent
\rom{(iii)}
The $\L(\pair)$
contains $S_{[\lm-\mu]}$ 
with multiplicity one 
as a $\C[W_\l]$-module.
\end{th}
\begin{rem}
  The statement (i) easily follows from (ii) and (iii).
\end{rem}
\begin{th}[\cite{Zel;induced}]
\label{th;classification}
Any simple $H_\l$-module of level $n$ is isomorphic to
$\L(\lm,\mu)$ for some
$\lm\in D_n$ and $\mu\in \lm-P(\Vn)$.
\end{th}
For $\lm\in D_n$ and $\mu\in\lm-P(\Vn)$,
 the $H_\l$-module $\K(\pair)$ is called
a {\it standard module}. 
For $\eta\in\h_n^*$, let $W_n[\eta]$ denote the stabilizer of $\eta$:
\begin{equation}
  \label{eq;stabilizer}
  W_n[\eta]=\{w\in W_n\mid w(\eta)=\eta\},
\end{equation}
which is a parabolic subgroup of $W_n$.
\begin{pr}[\cite{Zel;induced}]
\label{pr;isom_type}
Suppose that $\lm,\mu\in D_{n}$
and $w,y\in W_{n}$ satisfy
$\lm-w\circ\mu\in P(\Vn)$ and $\lm-y\circ\mu\in P(\Vn)$.
Then the following conditions are equivalent$:$

\medskip\noindent
\rom{(i)} $w\in W_n[\lm+\rho]yW_n[\mu+\rho]$.

\smallskip\noindent
\rom{(ii)} $\K(\lm,w\circ\mu)\cong\K(\lm,y\circ\mu)$.

\smallskip\noindent
\rom{(iii)} $\L(\lm,w\circ\mu)\cong \L(\lm,y\circ\mu)$.
\end{pr}
\begin{rem}\label{rem;important}
Let $\lm,\mu\in D_n$ and $w\in W_n$ such that
$\lm-w\circ\mu\in P(\Vn)$.
We often use the following fact from
Proposition \ref{pr;isom_type}:   
\begin{alignat}{3}
&\K(\lm,w\circ\mu)&\cong&\K(\lm,w^\lm\circ\mu)&\cong
&\K(\lm,w^\lm_\mu\circ\mu),\label{eq;isom_standard}\\
&\L(\lm,w\circ\mu)&\cong&\ \L(\lm,w^\lm\circ\mu)&\cong
&\L(\lm,w^\lm_\mu\circ\mu).\label{eq;isom_simple}
\end{alignat}
Here $w^\lm$ (resp. $w^\lm_\mu$) denotes the unique longest 
element in
$W_n[\lm+\rho]w$ (resp. $W_n[\lm+\rho]w W_n[\mu+\rho]$).
\end{rem}
\section{Functors $F_\lm$}
\subsection{Construction}
Let us recall the definition of the functor 
$$F_\lm:\O(\g_{n})\to\Rep(H_\l)$$
introduced in \cite{AS}.
Here $\Rep(H_\l)$ denotes the category of finite-dimensional representations of
$H_\l$.
Let $V_n=\C^n$ denote the vector
representation of $\g_n$. 
\begin{lem}[\cite{AS}]\label{lem;H-action}
For any $X\in\O(\g_n)$,
there exists a unique homomorphism
\begin{equation}
  \theta:H_\l \to\End_{\g_{n}}(X\*\Vn)
\end{equation}
such that
\begin{alignat}{2}
  s_i&\mapsto \Omega_{i\, i+1}&\quad &(i=1,\dots,\l-1),\\
  \e\ch_i&\mapsto \sum_{0\leq j<i}\Omega_{ji}+\frac{n-1}{2}&
\quad &(i=1,\dots,\l),\label{eq;act_S(t)}
\end{alignat}
where
\begin{align*}
\Omega_{ji}
           &=\sum_{1\leq k,m \leq n}
1^{\* j}\* E_{k,m}\* 1^{\* i-j-1}\*
E_{m,k}\* 1^{\* \l-i}
\in \End (X\*\Vn).
\end{align*}
\end{lem}
Let $\lm\in D_{n}$ and $X\in \ob \O(\g_n)$.
We define
\begin{equation}
  F_\lm(X)=(X\*\Vn)^{[\lm]}_\lm
\end{equation}
with an induced $H_\l$-module structure
through
the homomorphism $\theta$.
We also introduce an $H_\l$-module structure
on $\left((X\*\Vn)/\n^-_n(X\*\Vn)\right)_\lm$.
Then the bijection given in Lemma \ref{lem;bij}
gives an $H_\l$-isomorphism
\begin{equation}\label{eq;quo}
  F_\lm(X)\cong \left((X\*\Vn)/\n^-_n(X\*\Vn)\right)_\lm.
\end{equation}
Obviously $F_\lm$ defines an exact functor from $\O(\g_n)$ to 
$\Rep(H_\l)$.
\subsection{Image of functors}
We extend the definition of $\K(\pair)$
for any $\lm,\mu\in\t_n^*$ by
\begin{equation}\label{convention}
\K(\pair)=0 \text{ for }\lm,\mu\in\t_n^*\text{ such that }
\lm-\mu\notin P(\Vn).
\end{equation}
Let $\{u_i\}_{i=1,\dots,n}$ denote the standard basis of 
$V_n=\C^n$.
For $\lm\in D_n$ and $\mu\in\lm-P(\Vn)$, we
define an element
$u_{\lm,\mu}\in \left((M(\mu)\*\Vn)\left/\right.
\n_n^-(M(\mu)\*\Vn)\right)_\lm$
as the image of $v_{\mu}\* u_1^{\*\l_1}\*\cdots\* u_n^{\*\l_n}\in
M(\mu)\* \Vn$, where $\l_i=\bra\lm-\mu,\e\ch_i\ket_n$.
It was shown in \cite{AS} that
there exists an $H_\l$-homomorphism
\begin{equation}\label{eq;hom_standard}
  \K(\pair)\to \left(M(\mu)\*\Vn\left/\right.
\n_n^-(M(\mu)\*\Vn)\right)_\lm,
\end{equation}
which sends $\one_\pair$ to
$u_{\pair}$,
and that this is bijective.
Combining \eqref{eq;quo}, we have
\begin{th}[\cite{AS}]\label{th;image_of_Verma}
For each $\lm\in D_{n}$ and $\mu\in \h^*_{n}$,
there is an isomorphism of $H_{\l}$-modules
\begin{equation*}
  \label{eq;verma-standard}
  F_\lm(M(\mu))\cong \K(\pair).
\end{equation*}
In particular,
the $H_\l$-module $F_\lm(M(\mu))$ has the unique simple quotient. 
\end{th}
A proof of the
following statement is given  in \S \ref{ss;pf_simplicity}.
\begin{th}  \label{th;image_of_simple}
Let  $\lm\in D_n$ and 
$\mu\in \lm-P(\Vn)$.

\medskip\noindent \rom{(i)} If $\mu$ satisfies the condition
\begin{equation}
\label{eq;nonzero_condition}
  \bra \mu+\rho, \al\ch\ket_n\leq 0\quad \text{ for any } 
\al\in R^+_n
\text{ such that } \bra \lm+\rho,\al\ch\ket_n=0,
\end{equation}
then we have
\begin{equation}
   F_\lm(L(\mu))\cong \L(\lm,\mu),
\end{equation}
where $\L(\pair)$ is the unique simple quotient of $\K(\pair)$.

\smallskip\noindent
\rom{(ii)} If $\mu$ 
does not satisfy the condition \eqref{eq;nonzero_condition},
then we have
\begin{equation}
   F_\lm(L(\mu))=0.
\end{equation}
\end{th}
\begin{rem}\label{rem;equivalent}
\rom{(i)}
In the case $\l=n$, Theorem \ref{th;image_of_simple}
was proved in \cite{AS}
using the Kazhdan-Lusztig type multiplicity formula
for $\O(\g_n)$ and that for $\Rep(H_\l)$ (see \S \ref{ss;KL}).

\smallskip\noindent
\rom{(ii)} 
Recall that 
$W_n[\eta]\subseteq W_n$ denotes the stabilizer of $\eta\in\h_n^*$
(see \eqref{eq;stabilizer}).
Let $W_n^\eta$ denote the integral Weyl group of $\eta$:
\begin{equation}
  \label{eq;integral_Weyl_group}
  W_n^\eta=\{ w\in W_n\mid w\circ\eta-\eta\in Q_n\}.
\end{equation}
(Recall that $Q_n=\+_{i=1}^{n-1}\Z \al_i$.)
We can express $\mu$ in Theorem \ref{th;image_of_simple}
as 
$$
\mu=w\circ\tilde{\mu}$$
 with $\tilde{\mu}\in D_n$
and $w\in W_n^{\tilde{\mu}}$. 
Then the condition \eqref{eq;nonzero_condition}
is equivalent to 
\begin{equation*}
{\mu}=w^{\lm}\circ\tilde{\mu}
\quad\text{or equivalently}\quad 
\mu=w_{\tilde{\mu}}^{\lm}\circ\tilde{\mu}.
\end{equation*}
Here 
$w^{\lm}$ (resp. $w_{\tilde{\mu}}^{\lm}$)
denotes the unique
longest element in the coset 
$W_n[\lm+\rho]w$ (resp.
 $W_n[\lm+\rho]w W_n[{\tilde{\mu} +\rho}]$).
(Note that 
$w^\lm\circ\tilde\mu
=w_{\tilde{\mu}}^\lm\circ\tilde\mu$.)
\end{rem}

From Theorem \ref{th;classification} and
 Proposition \ref{pr;isom_type}, we have
\begin{cor}
  Any finite-dimensional 
simple $H_\l$-module of level $n$ is isomorphic to
$F_\lm(L(\mu))$ for some 
$\lm\in D_n$ and
$\mu\in \lm-P(\Vn)$ satisfying \eqref{eq;nonzero_condition}.
\end{cor}
\section{Contravariant forms and the Jantzen filtration}
\label{s;contra}
We remark on contravariant bilinear forms on
$\g_n$-modules and those on $H_\l$-modules.
We relate them via the functor $F_\lm$.
As a consequence,
we have a proof of
Theorem \ref{th;image_of_simple}
(a similar argument can be seen in 
the theory of Jantzen's translation 
functors \cite{Ja}).
We also prove
that the Jantzen filtration on
the Verma modules are transformed to
the Jantzen filtration on the standard modules.
\subsection{Contravariant forms}
Let $X\in\ob\O(\g_n)$.
A bilinear form $(~|~)_X:X\times X\to\C$ 
is called a $\g_n$-{\it contravariant form} if
\begin{equation}
  \label{eq;contra_sl}
   (xv|u)_X=(v|\sigma(x)u)_X\quad
 \text{\rom{for all }} 
u,v\in X,\ x\in\g_{n},
\end{equation}
where $\sigma$ is the transposition (\S \ref{ss;sl}).
For $Y\in\ob\Rep(H_\l)$,
a  bilinear form $(~|~)_Y:Y\times Y\to\C$ 
is called an $H_\l$-{\it contravariant form} if
\begin{equation}
  \label{eq;contra_H}
   (xv|u)_Y=(v|\invo(x)u)_Y\quad
 \text{\rom{for all }} 
u,v\in Y,\ x\in H_\l,
\end{equation}
where $\invo$ is given in Lemma \ref{lem;involution}.

Let us recall some fundamental facts on 
contravariant bilinear forms.
The following lemma is easily shown.
\begin{lem}\label{lem;properties_of_bilinear_form}
\rom{(i)}
Let $X\in\ob \O(\g_n)$ be equipped
with a $\g_n$-contravariant bilinear form $(~\mid ~)_X$. 
Then we have
\begin{align}
X^{[\lm]}\perp X^{[\mu]}\quad
&\text{\rom{ unless  }}\ \lm \in W_n\circ\mu,\label{eq;ortho1}\\
X_\lm\perp X_{\mu}\quad
&\text{\rom{ unless }}\ \lm =\mu.\label{eq;ortho2}
\end{align}

\smallskip\noindent
\rom{(ii)} Let $Y\in \ob\Rep(H_\l)$ 
be equipped with an $H_\l$-contravariant bilinear form $(~|~)_Y$.
Then we have
\begin{equation}\label{eq;ortho_for_Hecke}
  Y^{\gen}_\zeta\bot Y^{\gen}_{\eta}\ 
\text{\rom{ unless }}\ \zeta =\eta.
\end{equation}
\end{lem}
\begin{lem}
  \rom{(i)}
Let $\mu\in\h_n^*$. A $\g_n$-contravariant form on $M(\mu)$
is unique up to constant multiples.

\smallskip\noindent
\rom{(ii)}
Let $\lm\in D_n$ and $\mu\in \lm-P(\Vn)$.
 An $H_\l$-contravariant form on $\K(\pair)$
is unique up to constant multiples.
\end{lem}
\begin{pf}
(i) is well-known.
  We will prove (ii). 
Recall the decomposition \eqref{eq;decomposition_as_W-module}:
\begin{equation*}
  \K(\pair)\cong S_{[\lm-\mu]}\oplus
\bigoplus_{\beta\triangleright [\lm-\mu]}S_\beta^{\+ a_\beta}
\end{equation*}
as a $\C[W_\l]$-module.
Because an $H_\l$-contravariant form is 
$W_\l$-invariant,
its restriction 
to $S_{[\lm-\mu]}$ is unique up to constant, and we have
\begin{equation}\label{eq;bot}
  S_{[\lm-\mu]}\bot
   \+_{\beta\triangleright [\lm-\mu]}S_\beta^{\+ a_\beta}.
\end{equation}
From Theorem \ref{th;usq}-(i), $S_{[\lm-\mu]}$
generates $\K(\pair)$ over $H_\l$.
Thus the statement follows.
\end{pf}
It is easy to construct
a non-zero $\g_n$-contravariant form on $M(\mu)$. 
It is also known that there exists
a non-zero contravariant form
on $\K(\pair)$
(see \cite{Ro,CG} and also Remark \ref{rem;existence}).
In the rest of this paper, we
fix a canonical $\g_n$-contravariant form $(~|~)_{M(\mu)}$
on $M(\mu)$ 
by
$(v_\mu|v_\mu)_{M(\mu)}=1$.
The following lemma is easily shown.
\begin{lem}\label{lem;radical}
\rom{(i)} Let $\mu\in\h_n^*$ and
let $N$ be a unique maximal submodule of $M(\mu)$.
Then
\begin{equation}
  N=\text{\rom{rad}}(~\mid ~)_{M(\mu)},
\end{equation}
where {\rom{rad}}$(~\mid ~)_{M(\mu)}$ denotes the radical of 
$(~\mid ~)_{M(\mu)}$.

\medskip\noindent
\rom{(ii)}
Let $\lm\in D_n$ and $\mu\in\lm-P(\Vn)$.
Let $(~\mid ~)_{\K(\pair)}$ be a non-zero
$H_\l$-contravariant form
on $\K(\pair)$ and let 
$\cal{N}$ be a unique maximal submodule of $\K(\pair)$. 
Then we have 
$$\cal{N}=\text{\rom{rad}}(~\mid ~)_{\K(\pair)}.$$
\end{lem}
\begin{pf} 
(i) is well-known. Let us prove (ii).
It is obvious that rad$(~\mid ~)\subseteq \cal{N}$.
Theorem \ref{th;usq} implies that $\cal{N}\subseteq
 \+_{\beta\triangleright [\lm-\mu]}S_\beta^{\+ a_\beta}$
with some $a_\beta\in\Z_{\geq 0}$.
Thus we have
$S_{[\lm-\mu]}\bot \cal{N}$ by \eqref{eq;bot}.
Hence Theorem \ref{th;usq}-(i) implies that
$\cal{N}\subseteq$rad$(~|~)_{\K(\pair)}.$
\end{pf}
Let $X,Y\in\ob\O(\g_n)$ with $\g_n$-contravariant forms
$(~|~)_X,~ (~|~)_Y$. 
Then the tensor product
$X\* Y$ is equipped with a natural $\g_n$-contravariant bilinear form 
such that $(u\* v\mid u'\*v')_{X\* Y}=(u\mid u')_X\,(v\mid v')_Y$
for $u,u'\in X$ and $v,v'\in Y$.
The following simple lemma will play a key
role.
\begin{lem}\label{lem;contra_induce}
Let $\lm\in D_n$.
Let $X$ be a highest weight module 
$($i.e. a quotient of a Verma module$)$ of $\g_n$.

\medskip\noindent
\rom{(i)} The $\g_n$-contravariant form
on $X\* \Vn$ is also $H_\l$-contravariant, and thus it
induces an
$H_\l$-contravariant form on 
$(X\*\Vn)^{[\lm]}_\lm=F_\lm(X)$. 

\smallskip\noindent
\rom{(ii)} If the $\g_n$-contravariant form on $X$ is 
non-degenerate, then the induced
contravariant form on $F_\lm(X)$ is non-degenerate.
\end{lem}
\begin{pf}
(i) can be easily checked.
(ii) follows from Lemma \ref{lem;properties_of_bilinear_form}.
\end{pf}
As a consequence of Lemma \ref{lem;contra_induce}-(i), 
the canonical $\g_n$-contravariant form on $M(\mu)$ induces
an $H_\l$-contravariant form on 
$\K(\pair)=
F_\lm(M(\mu))$, which we call the canonical contravariant
form on $\K(\pair)$.
By Lemma \ref{lem;radical}-(i),
the $\g_n$-contravariant form on $L(\mu)$ is
 non-degenerate, and
it induces a non-degenerate 
 $H_\l$-contravariant form on
$F_\lm(L(\mu))$ by Lemma \ref{lem;contra_induce}-(ii). 
By Lemma \ref{lem;radical}-(ii), we have
\begin{cor}\label{cor;simple_or_zero}
Suppose that $\lm\in D_n$ and $\mu\in \lm-P(\Vn)$.
Then
the $H_\l$-module $F_\lm(L(\mu))$ is simple unless it is zero.
\end{cor}
\subsection{Proof of Theorem \ref{th;image_of_simple}}
\label{ss;pf_simplicity}
By $F_\lm(M(\mu))\cong \K(\pair)$ and Corollary 
\ref{cor;simple_or_zero},
it follows that $F_\lm(L(\mu))$ is isomorphic to $\L(\pair)$ 
or zero.
Hence the proof of Theorem \ref{th;image_of_simple} is reduced 
to the following lemma:
\begin{lem}\label{lem;nonvanishing}
Let $\lm\in D_n$  and $\mu\in \lm-P(\Vn)$.
Then $ F_\lm(L(\mu))\neq 0$
if and only if $\mu$ satisfies the condition \eqref{eq;nonzero_condition}.
\end{lem}
\begin{rem}\label{rem;existence}
  Lemma \ref{lem;nonvanishing} implies that
the canonical $\g_n$-contravariant form on $M(\mu)$
induces a non-zero $H_\l$-contravariant form on
$F_\lm(M(\mu))$ if and only if the condition 
\eqref{eq;nonzero_condition} is satisfied.
By Remark \ref{rem;important} and Remark \ref{rem;equivalent},
it follows that
any standard module admits a non-zero $H_\l$-contravariant form.
\end{rem}

\noindent
{\it Proof of Lemma \ref{lem;nonvanishing}.}
First we show the ``only if'' part.
Suppose that $\mu$ does not satisfy \eqref{eq;nonzero_condition}. 
Then there exists $\al\in R_n^+$ such that
$\bra \mu+\rho,\al\ch\ket\in \Z_{>0}$ and 
$\bra\lm+\rho,\al\ch\ket=0.$
The first inequality implies
$M(s_\al\circ\mu)\subset M(\mu)$, and the second equality
implies
$\K(\lm,\mu)\cong\K(\lm,s_\al\circ\mu)$ 
(Proposition \ref{pr;isom_type}).
Hence we have $F_\lm(L(\mu))=0$, because it is a quotient of
$F_\lm(M(\mu))/F_\lm(M(s_\al\circ\mu))=0$.

Let us prove the ``if'' part.
We can write $\mu$ as
$$\mu=w\circ\tilde{\mu},$$
where $\tilde{\mu}\in D_n$
and $w$ is an element of
the integral Weyl group $W_n^{\tilde{\mu}}$ 
(see \eqref{eq;integral_Weyl_group}).

Then the condition 
\eqref{eq;nonzero_condition} implies 
$\mu=w_{\tilde{\mu}}^\lm\circ\tilde{\mu}$,
where $w_{\tilde{\mu}}^\lm$ is 
the longest element in 
$W_n[\lm+\rho]w_{\tilde{\mu}}W_n[\tilde{\mu}+\rho]$
(see Remark \ref{rem;equivalent}).
In the Grothendieck group of $\O(\g_n)$, we write 
\begin{equation}\label{eq;decomp}
  M(w_{\tilde{\mu}}^{\lm}\circ\tilde{\mu})=L(w_{\tilde{\mu}}^{\lm}
\circ\tilde{\mu})
+\sum_{y_{\tilde{\mu}}}a_{y_{\tilde{\mu}}} L(y_{\tilde{\mu}}
\circ\tilde{\mu}).
\end{equation}
Here the sum runs over 
those elements $y_{\tilde{\mu}}\in W_n$ 
such that $y_{\tilde{\mu}}$ is longest in
$y_{\tilde{\mu}}W_n[\mu+\rho]$ and  
\begin{equation}
  \label{eq;ine}
 y_{\tilde{\mu}}>w_{\tilde{\mu}}^{\lm}.
\end{equation}
Applying $F_\lm$ to \eqref{eq;decomp}
we have 
\begin{equation}
  \K(\lm,w_{\tilde{\mu}}^{\lm}\circ\tilde{\mu})
=F_\lm(L(w_{\tilde{\mu}}^{\lm}\circ\tilde{\mu}))+
\sum_{y_{\tilde{\mu}}}a_{y_{\tilde{\mu}}} 
F_\lm(L(y_{\tilde{\mu}}\circ\tilde{\mu}))
\end{equation}
in the Grothendieck group of $\Rep(H_\l)$.
Assuming
 that $F_\lm(L(w_{\tilde{\mu}}^{\lm}\circ\tilde{\mu}))=0$,
we will deduce a contradiction.
Since the multiplicity of 
$\L(\lm,w_{\tilde{\mu}}^{\lm}\circ\tilde{\mu})$
in $\K(\lm,w_{\tilde{\mu}}^{\lm}\circ\tilde{\mu})$ 
is non-zero,
Corollary \ref{cor;simple_or_zero} implies
$$\L(\lm,w_{\tilde{\mu}}^{\lm}\circ\tilde{\mu})=
F_\lm(L(y_{\tilde{\mu}}\circ\tilde{\mu}))=
\L(\lm,y_{\tilde{\mu}}\circ\tilde{\mu})$$
for some $y_{\tilde{\mu}}$.
But this implies 
$y_{\tilde{\mu}}\in 
W_n[\lm+\rho]w_{\tilde{\mu}}^{\lm}W_n[\tilde\mu+\rho]$
by Proposition \ref{pr;isom_type},
and thus we have $l(y_{\tilde{\mu}})
\leq l(w_{\tilde{\mu}}^{\lm})$.
This contradicts \eqref{eq;ine}.\qed
\subsection{The Jantzen filtrations}\label{ss;Jantzen}
Throughout
 this subsection, we fix a weight $\delta\in\h^*_n$.
Let $A=\C[\,t \,]_{(t)}$ denote the localization of
$\C[\,t\,]$ at the prime ideal $(t)$. 
We use the notation: $\eta^t=\eta+\delta t\in\h^*_n\*A$
for $\eta\in \h^*_n$.

For $\mu\in \h_n^*$,
let $M(\mu^t)$ be the Verma module of $\g_n\* A$ 
with highest weight $\mu^t$:
\begin{align*}
&M (\mu^t)= (U(\g_n)\* A)\*_{U(\b^+_n)\* A}(Av_{\mu^t}).
\end{align*}
The canonical $\g_n$-contravariant bilinear form on $M(\mu)$ 
can be naturally
extended to 
a $\g_n\*A$-contravariant form $(~|~)_{M(\mu^t)}$
on $M(\mu^t)$ 
(with respect to the anti-involution
$\sigma\*\text{\rom{id}}_A$)
with values in $A$.

Define
\begin{align}
  M(\mu^t)_j&=\{v\in M(\mu^t)\mid (v\mid u)_{M(\mu^t)}\in 
t^j A\text{ for all } u\in M(\mu^t)\}.
 \end{align}
Putting  $M(\mu)_j=M(\mu^t)_j\left/\left( tM(\mu^t)\cap M(\mu^t)_j
\right)\right.$
we have a filtration
\begin{equation}
  M(\mu)=M(\mu)_0\supseteq M(\mu)_1\supseteq M(\mu)_2\supseteq
\cdots
\end{equation}
by $\g_n$-modules
called the {\it Jantzen filtration} \cite{Ja}. 

Our next aim is to define 
the Jantzen filtration on the standard module, which was 
introduced in \cite{Ro}.
Let $\lm\in D_n$ and $\mu\in\lm-P(\Vn)$.
Analogously to \S\ref{ss;induced}, we define
an $H_\l\* A$-module $\K(\lm^t,\mu^t)$ by
$$
\K(\lm^t,\mu^t)=(H_\l\* A)\*_{H_{\lm,\mu}\* A}
(A\one_{\lm^t,\mu^t}).
$$

Put $X=M(\mu^t)\* \Vn$, which is
 equipped with a $\g_n\* A$-contravariant 
form $(~|~)_X$.
Then $\h_n^*\* A$ acts semisimply on $X$ and it follows that
\begin{equation}
  X=\+_{\eta^t\in \mu^t+P_n}X_{\eta^t},
\end{equation}
\begin{equation}
  \label{eq;ortho1A}
  X_{\eta^t}\bot X_{\nu^t} \text{ unless }
\mu=\nu.
\end{equation}
Let $\chi_{\eta^t}:Z(U(\g_n)\* A)\to A$
be the infinitesimal character of $M(\eta^t)$.
Following \cite{GJ}, we
define for  $\eta\in\h_n^*$
an ideal $J_{\eta^t}$ of $Z(U(\g_n)\* A)$
by
 $$J_{\eta^t}=\cap_{w\in W_n}
 \Ker \chi_{(w\circ\eta)^t},$$
and define
\begin{equation}
  X^{[\eta^t]}=\{v\in X\mid J_{\eta^t}^kv=0\ \text{ for some }k\}.
\end{equation}
Obviously $X^{[\eta^t]}$ depends only on
the equivalence class $[\eta]$ of $\eta$ with
respect to the equivalence relation \eqref{eq;equivalence_relation}.
\begin{lem}[\cite{GJ}]\label{lem;decomposition}
We have
  \begin{equation}
    X=\+_{[\eta]\in \t_n^*/\sim}X^{[\eta^t]},
  \end{equation}
\begin{equation}
\label{eq;ortho2A}
  X^{[\eta^t]}\bot X^{[\nu^t]}\ \text{\rom{ unless }}
[\eta]=[\nu].
\end{equation}
\end{lem}

On the $\g_n\* A$-module
$X=M(\mu^t)\* \Vn$,
we can define an action of $H_\l\*A$ commuting with
$\g_n\* A$ as in Lemma \ref{lem;H-action}.
We define an induced
$H_\l\* A$-module structure on the following spaces:
\begin{align}
&(X/\n_n^-X)_{\lm^t},\quad
(X^{[\lm^t]})_{\lm^t}.
\end{align}
With respect to this action,
the natural map
\begin{equation}\label{eq;hom1}
(X ^{[\lm^t]})_{\lm^t}\to (X/\n_n^- X)_{\lm^t}
\end{equation}
is an $H_\l\* A$-homomorphism.

Similarly to \eqref{eq;hom_standard}, 
we can construct an $H_\l\* A$-homomorphism
\begin{equation}
  \label{eq;hom2}
  \K(\lm^t,\mu^t)\to (X/\n_n^- X)_{\lm^t}.
\end{equation}
The following lemma is elementary.
\begin{lem}\label{lem;Nakayama}
  Let $M$ and $N$ be free $A$-modules of finite rank,
and let 
$f:M\to N$
be an $A$-homomorphism.
If the specialization 
$$\bar{f}:M/tM\to N/tN$$
at $t=0$
is a $\C$-isomorphism, then $f$ is an $A$-isomorphism.
\end{lem}
Using Lemma \ref{lem;Nakayama}, we get
\begin{pr}
  The $H_\l\* A$-homomorphisms \eqref{eq;hom1} and \eqref{eq;hom2}
are bijective$:$
\begin{equation}
(X^{[\lm^t]})_{\lm^t}
\cong (X/\n_n^- X)_{\lm^t}
\cong  \K(\lm^t,\mu^t).
\end{equation}
\end{pr}
\begin{pf}
The specialization of  \eqref{eq;hom1} (resp. \eqref{eq;hom2})
at $t=0$ gives the isomorphism in 
Lemma \ref{lem;bij} (resp. \eqref{eq;hom_standard}).
Therefore by 
Lemma \ref{lem;Nakayama}, 
it is enough to show that 
$(X^{[\lm^t]})_{\lm^t}$, 
$(X/\n_n^- X)_{\lm^t}$
and $\K(\lm^t,\mu^t)$ are
all free $A$-modules of finite rank.
Obviously they are finitely generated over $A$.
 It is also clear that $\K(\lm^t,\mu^t)$ is free.
Since $A$ is a principal ideal domain
and $X$
is a free $A$-module,
its subspace $(X^{[\lm^t]})_{\lm^t}$
is a free $A$-module.
Finally, let us show that
$(X/\n_n^- X)_{\lm^t}$
is a free $A$-module.
By the isomorphism 
\begin{equation}
X= M(\mu^t)\* \Vn\cong (U(\g_n)\* A)\*_{U(\b_n^+)\* A}
(A v_{\mu^t}\* \Vn)
\end{equation}
as $U(\g_n)\* A$-modules, it follows that
\begin{equation}
 (X/\n_n^- X)_{\lm^t}\cong (\Vn)_{\lm-\mu}\* A  
\end{equation}
as $A$-modules. This is a free $A$-module.
\end{pf}

It follows that the $\g_n\* A$-contravariant form on 
$X=M(\mu^t)\*\Vn$
is also $H_\l\*A$-contravariant.
Through the isomorphism
\begin{equation}
  \K(\lm^t,\mu^t)\cong 
(X^{[\lm^t]})_{\lm^t}\subset X,
\end{equation}
we introduce an $A$-valued 
$H_\l\* A$-contravariant form on $\K(\lm^t,\mu^t)$.

Assume that
 $\mu$ satisfies the condition \eqref{eq;nonzero_condition}
in Theorem \ref{th;image_of_simple}. Then
the induced contravariant form is non-zero (since
its specialization at $t=0$ is non-zero).
Therefore we have a filtration
\begin{equation}
\K(\pair)=  \K(\pair)_0\supseteq \K(\pair)_1\supseteq
\K(\pair)_2\supseteq\cdots
\end{equation}
by $H_\l$-modules, which we call the Jantzen filtration.
Recall that any standard module is isomorphic to
$\K(\lm,\mu)$ for some $\lm\in D_n$ and $\mu\in\lm-P(\Vn)$
satisfying \eqref{eq;nonzero_condition}
(Remark \ref{rem;important}).
\begin{rem}
In \cite{Ro}, the deformation direction $\delta$ is restricted by
a certain condition. The construction above gives
the definition of the Jantzen filtration for an arbitrary 
direction $\delta$.
\end{rem}
\begin{th}\label{th;Jantzen}
Suppose that $\lm\in D_n$ and $\mu\in\lm-P(\Vn)$
satisfy the condition
\eqref{eq;nonzero_condition}.
Then $F_\lm(M(\mu)_j)=\K(\pair)_j$.
\end{th}
\begin{pf}
It is easy to check that
$F_\lm(M(\mu)_j)\subseteq \K(\pair)_j$.
To prove the opposite inclusion, let 
$$p:M(\mu^t)\* \Vn\to  (M(\mu^t)\*\Vn)^{[\lm^t]}_{\lm^t}
=\K(\lm^t,\mu^t)$$
denote the natural projection. 
Note that $ (M(\mu^t)\*\Vn))^{[\lm^t]}_{\lm^t}\bot\Ker p$
by  \eqref{eq;ortho1A} and \eqref{eq;ortho2A}.
Fix any orthonormal basis $\{b_i\}_{i=1}^{n^\l}$ of
$\Vn$ with respect to the $\g_n$-contravariant form 
$(~\mid ~)_{\Vn}$.

Take any $u\in \K(\lm^t,\mu^t)_j\subseteq
(M(\mu^t)\*\Vn))^{[\lm^t]}_{\lm^t}$
and write as $u=\sum_i a_i\* b_i$ with $a_i\in M (\mu^t)$.
Then for any $v\in M(\mu^t)$ and $k$, 
we have
\begin{align*}
(a_k\mid v)_{M(\mu^t)}&=(u\mid v\* b_k)_{M(\mu^t)\*\Vn}=
(u\mid p(v\* b_k))_{M(\mu^t)\*\Vn}\\
&=(u\mid p(v\* b_k))_{ (M(\mu^t)\*\Vn)^{[\lm^t]}_{\lm^t}}
\in t^j A.
\end{align*}
This implies $a_k\in M(\mu^t)_j$ and thus
$u\in (M(\mu^t)_j\*\Vn)^{[\lm^t]}_{\lm^t}.$  
Therefore we have $F_\lm(M(\mu)_j)\supseteq \K(\pair)_j$.
\end{pf}
\section{Consequences}\label{s;consequences}
\subsection{BGG resolution}
Recall the generalization of the BGG resolution 
for certain simple $\g_n$-modules given by
Gabber-Joseph \cite{GJ;BGG}.

We fix $\mu\in \h^*_n$ such that
$-(\mu+\rho)$ is dominant and regular,
i.e. $\bra -(\mu+\rho),\al\ch\ket_n\notin \Z_{\leq 0}$
for all $\al\in R_n^+$.
Set $R^\mu_n=\{\al\in R_n\mid \bra\mu,\al\ch\ket_n\in \Z\}$.
It is known that $R_n^\mu$ is a root system and its
Weyl group coincides with the integral Weyl group
\begin{equation}
 \label{eq;integral_Weyl_group2}
W^\mu_n=\{w\in W_n\mid w\circ\mu-\mu\in Q_n\}.
  \end{equation}
Set $R_n^{\mu+}=R_n^\mu\cap R_n^+$ and 
let $\Pi_n^\mu$ be the
set of simple roots of $R_n^{\mu+}$.

Fix $B \subseteq \Pi_n^\mu$.
The length function $l_B$ and the Bruhat order
of $W_B$ are defined
with respect to the set of simple roots $B$.
Let $w_B$ be a unique longest element of $W_B$
with respect to $l_B$.
Put $\mu_B=w_B\circ\mu$.
Gabber-Joseph constructed the exact sequence
\begin{equation}
  \label{eq;BGG2}
0\gets  L(\mu_B)\gets C_0\gets C_1\gets\cdots
\end{equation}
of $\g_n$-modules, where
$$C_i=\+
\begin{Sb}
y\in W_B,\
l_B(y)=i
\end{Sb}
M(y\circ\mu_B).$$
We apply $F_\lm$ to the sequence \eqref{eq;BGG2}.
Then Theorem \ref{th;image_of_Verma}
and Theorem \ref{th;image_of_simple}
imply the following:
\begin{th}\label{th;BGG_for_Hecke}
  Let $\mu$ and $B$ as above.
Suppose that $\lm\in D_n\cap (\mu_B+P(\Vn))$ satisfies
$\bra \lm+\rho,\al\ch\ket\neq 0$ for 
any $\al\in B$.
Then 
there exists an exact sequence 
\begin{equation}
0\gets \L(\lm,\mu_B)\gets 
{\cal C}_0\gets {\cal C}_1\gets \cdots
\end{equation}
of $H_\l$-modules, where 
$${\cal C}_i=\+
\begin{Sb}
y\in W_B,\ l_B(y)=i
\end{Sb}
\K(\lm,y\circ\mu_B).$$
\end{th}
\begin{rem}
In the case $\mu_B\in P_n^+$ and
$B=\Pi_\l$ (the original BGG case \cite{BGG}), 
the corresponding sequence has been obtained by Cherednik 
\cite{Ch;character} by a different method
(see also \cite{Zel;resolvant,AST}).
\end{rem}
\subsection{Kazhdan-Lusztig formulas}\label{ss;KL}
For a module $M$ and simple module $L$,
let $[M:L]$ denote the multiplicity of $L$
in the composition series of $M$.

Recall that $W_n^\mu$ denotes the integral Weyl group of $\mu\in\h^*_n$
(see  \eqref{eq;integral_Weyl_group}).
The following formula is
a direct consequence of Theorem \ref{th;image_of_Verma} and 
Theorem \ref{th;image_of_simple}:
\begin{th}\label{th;multiplicity}
Let $\lm, \mu\in D_n$
and  let $w, y\in W_n^\mu$
such that $\lm-w\circ\mu,\lm-y\circ \mu\in P(\Vn)$.
Then we have
  \begin{equation}\label{eq;multiplicity}
    [\K(\lm,w\circ \mu):\L(\lm,y\circ \mu)]=
    [M(w\circ\mu):L(y^{\lm}\circ\mu)],
  \end{equation}
where
$y^{\lm}$ denotes the longest element in
$W_n[\lm+\rho]y$. 
\end{th}
Let $\lm,\mu\in D_n$ and $w,y\in W_n^\mu$ be as in 
Theorem \ref{th;multiplicity}.
The equality \eqref{eq;multiplicity}
has been known 
through the following two multiplicity formulas:
\begin{alignat}{2}\label{eq;KL}
  [M(w\circ \mu)&:L(y\circ \mu)]& &=P_{w, y_{\mu} }(1),\\
\label{eq;KL_for_Hecke}
   [\K(\lm,w\circ\mu)&:\L(\lm,y\circ\mu)]&
&=P_{w, y_{{\mu}}^{\lm}}(1).
\end{alignat}
Here $P_{w,y}(q)\in \Z[q,q^{-1}]$ denotes the 
Kazhdan-Lusztig polynomial \cite{KL;Coxeter groups}
of the Hecke algebra associated to $W_n^\mu$ 
(we put $P_{w,y}(q)=0$ for $w\not< y$ for convenience), and
 $y_{\mu}$ (resp. $y_{\mu}^{\lm}$) denotes 
the longest element in $yW_n[\mu+\rho]$
(resp. $W_n[\lm+\rho]yW_n[\mu+\rho]$).
\begin{rem}
It follows from \eqref{eq;KL} and \eqref{eq;KL_for_Hecke}
that
$P_{w,y_\mu}(1)=P_{w_\mu,y_\mu}(1)$ and
$P_{w,y_{\mu}^\lm}(1)=P_{w_\mu,y_\mu^\lm}(1)
=P_{w_\mu^\lm,y_\mu^\lm}(1)$.
The latter is expressed
in terms of the intersection cohomology concerning 
nilpotent orbits on the quiver variety \cite{Zel;two remarks}.
\end{rem}

The formula \eqref{eq;KL} was conjectured by 
Kazhdan-Lusztig\ \cite{KL;Coxeter groups}
and proved by Beilinson-Bernstein\ \cite{BB}
and Brylinski-Kashiwara\ \cite{BK}.
The formula \eqref{eq;KL_for_Hecke} 
was conjectured by Zelevinsky\ \cite{Zel;p-adic KL}
(see also \cite{Zel;two remarks})
and proved by Ginzburg\ \cite{Gi} (see also \cite{CG}).
The theory of perverse sheaves plays an essential role
in these proofs.

Theorem \ref{th;multiplicity} (proved in a purely algebraic way)
says that the
Kazhdan-Lusztig formula \eqref{eq;KL}
is equivalent to its degenerate affine Hecke analogue 
(or its p-adic analogue)
\eqref{eq;KL_for_Hecke}.
The implication $\eqref{eq;KL}\Rightarrow \eqref{eq;KL_for_Hecke}$
is obvious.
The implication $\eqref{eq;KL_for_Hecke}
\Rightarrow  \eqref{eq;KL}$ is proved as follows.
Take any $\mu\in D_n$ and $w,y\in W_n^\mu$.
Then we can find $\l\in \Z_{\geq 2}$ and $\lm\in D_n^\circ$
such that 
$$  \lm-z\circ\mu\in P(\Vn)
\text{ for all }z\in W_n^\mu.$$
In this case
$F_\lm(L(z\circ\mu))$ never vanishes
and thus it is isomorphic to $\L(\lm,z\circ\mu)$.
Now \eqref{eq;KL_for_Hecke}
implies \eqref{eq;KL}.
\subsection{Rogawski's conjecture}\label{ss;Rog}
Let $\{M(\mu)_j\}_j$ and $\{\K(\pair)_j\}_j$
be the Jantzen filtrations defined in
\S \ref{ss;Jantzen}.
As a direct consequence of 
Theorem \ref{th;image_of_simple}
and Theorem \ref{th;Jantzen}, we have
\begin{th}\label{th;improved_multiplicity}
Let $\lm,\mu\in D_n$ and $w,y\in W_n^\mu$ 
$($see \eqref{eq;integral_Weyl_group}$)$ be 
such that $\lm-w\circ\mu,\lm-y\circ \mu\in P(\Vn)$.
Then we have
  \begin{equation}
    [\K(\lm,w\circ \mu)_j:\L(\lm,y\circ \mu)]=
    [M(w^{\lm}\circ\mu)_j:L(y^{\lm}\circ\mu)],
  \end{equation}
where $w^\lm$ and $y^{\lm}$ denote the longest element in
$W_n[\lm+\rho]w$ and $W_n[\lm+\rho]y$ respectively. 
\end{th}
A priori the Jantzen filtrations depend on the choice
of the deformation direction $\delta\in\h_n^*$.
It has been known
that
the Jantzen filtration on $M(\mu)$
does not depend on the choice of
$\delta$ for which $(~|~)_{M(\mu^t)}$ is non-degenerate
\cite{Ba}. 
Now Theorem \ref{th;Jantzen} implies  
\begin{pr}\label{pr;indep}
Let $\lm\in D_n$ and $\mu\in \lm-P(\Vn)$ satisfy 
\eqref{eq;nonzero_condition}.
Then the Jantzen filtration on $\K(\pair)$ does not depend on
the choice of $\delta$ such that
\begin{equation}
  \label{eq;NDC}
  \bra\delta,\al\ch\ket_n\neq 0\ \text{{ for any }}
\al\in R_n^+ \text{{ such that }}
  \bra \mu+\rho,\al\ch\ket_n\in \Z_{> 0}.
\end{equation}
\end{pr}
\begin{rem}
For $\lm$ and $\mu$ as in Proposition \ref{pr;indep},
  the condition \eqref{eq;NDC} is equivalent to
the condition that the $H_\l\* A$-contravariant form
$(~|~)_{\K(\lm^t,\mu^t)}$ is non-degenerate.
\end{rem}
We say that the Jantzen filtration $\{M(\mu)_j\}_j$
(or $\{\K(\pair)_j\}_j$) is {\it regular} if
the deformation direction $\delta$ satisfies \eqref{eq;NDC}.
The following formula
was conjectured in \cite{GJ,GM}, and
 proved in \cite{BB;Jantzen}.
\begin{th}[\cite{BB;Jantzen}]\label{th;BB}
Let $\mu\in D_n$ and $w,y\in W_n^\mu$. 
Suppose that $w$ and $y$ are the longest elements in 
$wW_n[\mu+\rho]$ and $yW_n[\mu+\rho]$, respectively.
For the regular Jantzen filtration $\{M(w\circ\mu)_j\}_j$,
we have
  \begin{equation}\label{eq;BB}
    \sum_{j\in \Z_{\geq 0}} 
[\gr_j M(w\circ\mu):L(y\circ\mu)]q^{(l_\mu(y)-l_\mu(w)-j)/2}
=P_{w,y}(q),
  \end{equation}
where $P_{w,y}(q)$ denotes the Kazhdan-Lusztig polynomial
of $W_n^\mu$,
and $l_\mu$ denotes the length function on
$W_n^\mu$.
\end{th}
Combining with Theorem \ref{th;improved_multiplicity},
the improved Kazhdan-Lusztig formula \eqref{eq;BB}
implies its degenerate affine Hecke analogue, 
which was conjectured in \cite{Ro}.
\begin{th}
${\hbox{\rom{\bf (}}}{\hbox{\rom{c.f.}}}$ 
\cite[Theorem 2.6.1]{Gi;ICM}${\hbox{\rom{\bf )}}}$
\label{th;Rog}
Let $\lm,\mu\in D_n$ and $w,y\in W_n^\mu$
be such that $\lm-w\circ\mu,~\lm-y\circ\mu\in P(\Vn)$.
Suppose that 
$w$ and $y$ are the longest elements in 
$W_n[\lm+\rho]wW_n[\mu+\rho]$ and $W_n[\lm+\rho]yW_n[\mu+\rho]$, 
respectively.
For the regular Jantzen filtration
 $\{\K(\lm,w\circ\mu)_j\}_j$,
we have
\begin{equation} 
  \sum_{j\in \Z_{\geq 0}} [\gr_j\K(\lm,w\circ\mu):\L(\lm,y\circ\mu)]
q^{(l_\mu(y)-l_\mu(w)-j)/2}
=P_{w,y}(q),
\end{equation}
where $P_{w,y}(q)$ denotes the Kazhdan-Lusztig polynomial
of $W_n^\mu$, and
$l_\mu$ denotes the length function on $W_n^\mu$. 
\end{th}

\end{document}